\documentclass[12pt]{article}

\usepackage{graphicx}
\usepackage{bm}

\usepackage[english]{babel}
\usepackage{amsmath,amssymb,amsthm,amsfonts}
\usepackage{mathrsfs}

\renewcommand{\geq}{\geqslant}
\renewcommand{\leq}{\leqslant}

\begin{document}

%\pagestyle{empty}

%  Include here your own macros.

\newcommand{\Eqref}[1]{(\ref{#1})}

%  Data for the headings. Please fill both fields.

%  Full title. Use '\\' to force a line break.

\title{Ecoepidemics with infected prey in herd defense: the harmless and toxic cases.}
%Ecoepidemics with group defense: the case of infected not segregated prey.}

%  Full author information.

%\author{Full name}{e-mail}{affiliation number (below)}

\author{Elena Cagliero, Ezio Venturino\\
%{\footnote {This paper was completed and written during a visit of the second author
%at the Max Planck Institut f\"ur Physik Komplexer Systeme in Dresden, Germany. The author expresses his thanks for the facilities provided.}
Dipartimento di Matematica ``Giuseppe Peano'',\\
Universit\`a di Torino, \\via Carlo Alberto 10, 10123 Torino, Italy,\\
emails: elenacagliero87@yahoo.it, ezio.venturino@unito.it}

\maketitle

\begin{abstract}
We consider a predator-prey population model with prey gathering together for defense purposes.
A transmissible unrecoverable disease affects the prey. We characterize the system behavior, establishing that
ultimately either only the susceptible prey survive, or the disease becomes endemic, but the predators are wiped out. Another
alternative is that the disease is eradicated, with sound prey and predators thriving at an equilibrium or through persistent
population oscillations. Finally, the populations can thrive together, with the endemic disease.
The only impossible alternative in these circumstances is predators thriving just with infected prey.
But this follows from the model assumptions, in that infected prey are too weak to sustain themselves.
A mathematical peculiarity of the model is the singularity-free reformulation, which leads to three entirely new dependent variables to describe
the system. The model is then extended to encompass the situation in which ingestion of diseased prey is fatal for the predators.

\end{abstract}

{\textbf{Keywords}: group defense, epidemics, predator-prey, disease transmission, toxic prey\\
{\textbf{AMS MR classification} 92D30, 92D25, 92D40
%
%  Main text of the article.
%

\section{Background}
%2 MODELLI, LE PREDE INFETTE SONO HARMLESS, OPPURE DANNEGGIANO I PREDATORI

The model we consider here is a prey-predator system in which the disease develops in prey. The latter gather together and live in a herd.
Following recently introduced ideas, \cite{AV,APV}, the large predators will hunt alone the herd and in it, it will be the individuals on
the edge of the bunch that will mostly bear the burden of the attack. In mathematical terms, the ``size'' of the prey population occupying
the edge of the herd is proportional to the square root of the total population. Thus, instead of the standard mass action or Holling type II
terms usually employed to model the predation mechanism, the predator-prey interactions are mathematically described via a term containing
the square root of the prey population, coupled as usual with the predators' population. This is a different idea
from the approach as the one used in \cite{FW},
in which the defense mechanism is modelled via a suitable response function.
In \cite{EV11}, these ideas are extended to another situation, in which a disease affects the prey. For further developments, see \cite{B}.
Thus, in this way the first ecoepidemic model of this sort is proposed. An idea of this kind had been presented for predators hunting in
packs in \cite{C-DA}.

Ecoepidemic models in fact contain a basic interacting population system on top of which a contagious disease is present.
Models of this type are known since about a quarter of a century, \cite{HF,BC,V95} and are currently of wide interest among scientists.
For an account of some of the developments of this branch of mathematical biology merging the two fields of population
theory and epidemiology, see \cite{h} or \cite{MPV}.

Coupling ecoepidemic systems with group defense is a very recent step, \cite{EV11}. In the formulation of the
model however, there is a kind of asymmetry in the way in which healthy and infected prey are dealt with by predators. Although both
are hunted, the predation assumes in \cite{EV11} two different mathematical forms, one containing the square root as discussed above,
the other one the standard Holling type I interaction term.
In fact,
the additional basic assumption with respect to the standard predator-prey model of \cite{APV}, which we will remove here, that has been
formulated in \cite{EV11} consists in the fact that the diseased prey are assumed to be left behind by the
healthy herd. Therefore they are subject to hunting by predators on a one-to-one basis, a
fact which is modeled as in the classical Lotka-Volterra system with the standard mass action term.

Here we want to extend the system studied in \cite{CV},
to encompass the situation in which the infected prey still remain in the herd, and mix with the
healthy ones. Therefore they can occupy any position in the bunch, including the ones near the boundary. They are therefore subject to hunt as all the
other susceptible prey. Mathematically speaking, the change amounts to the following:
the square root term that formerly contained only healthy individuals, is now replaced by a square root
term containing the whole prey population.

In this paper, we also extend the rephrased model in another direction. Instead
of considering the infected prey as a source of food for the predators as the healthy prey are, we also include the case
that the consumption of diseased prey has poisonous effects on the predators.

The paper is organized as follows. In the next Section we present the basic model.
In Section \ref{sec:ref} we redefine the original variables to obtain a singularity-free system and adimensionalize it.
The system's equilibria are assessed in Section \ref{sec:equil}.
Section \ref{sec:stab} contains their stability analysis. Hopf bifurcations are investigated in
Section \ref{sec:Hopf}.
% The brief 
Section \ref{sec:interp} summarizes the results interpretation in
terms of the original model variables. 
Section \ref{pois} deals with the case of predators being lethally affected by the consumption of diseased prey.
%Simulations are presented in Section \ref{sec:simul} and
A final discussion concludes the paper.

%\newpage
\section{The model}
Let $S$ denote the healthy prey population, $I$ be the infected prey and $P$ the predators.
We assume that the infection process running among the prey does not hinder them, so that infected individuals can still remain in the herd.
The predators attack the prey, and the individuals at the edge of the bunch are the most likely to be captured by the predators. Since the infected
do not remain behind the herd, they populate both the ``inside'' of the bunch as well as its boundary. Therefore they can be captured as well
as the healthy prey.

Following the arguments expounded in \cite{AV,APV,EV11}, if we
assume that the total prey population density $S+I$ is uniformly distributed on the land occupied by the herd, the number
of the individuals staying on the border is
proportional to the square root of this density.
With these assumptions the system can be written as
\begin{eqnarray}\label{sistema0}
%\begin{cases}
%\displaystyle{
\frac {dS}{dt} &=& rS\left( 1-\frac{S+I}{K} \right)-\sigma \frac{SI}{S+I}-qPS \frac{\sqrt{S+I}}{S+I} %}
\\ \nonumber
%\bigskip
%\displaystyle{
\frac {dI}{dt} &=& \sigma \frac{SI}{S+I}-w PI \frac{\sqrt{S+I}}{S+I}- \mu I %}
\\ \nonumber
%\displaystyle{
\frac {dP}{dt} &=& -mP+gP\frac{S}{\sqrt{S+I}} +fP\frac{I}{\sqrt{S+I}} 
%-mP +eP\sqrt{S+I}, %}
%\end{cases}
\end{eqnarray}
where all the parameters are assumed to be nonnegative.
Here, $r$ denotes the birth rate of healthy prey,
$\sigma$ is the disease incidence,
$q$ the predation rate on healthy prey,
$w$ the predation rate on infected prey,
$\mu$ the natural plus disease-related mortality rate of infected prey,
$m$ the death rate of predators,
$e$ is the uptake due to predation for the predators,
$K$ the environment's carrying capacity.

The first equation shows that healthy prey follow a logistic growth, with intraspecific competition due also to the infected. Then there is the
disease contagion mechanism, which is here assumed to be modelled by the standard incidence. Finally, healthy prey on the edge of the herd
are captured by the predators, at rate $q$. Note that the last term expresses how many sound prey stay on the border.
In fact, the population on the boundary is $\sqrt{S+I}$ as argued earlier. Of this, only the fraction $S(S+I)^{-1}$ is represented by healthy individuals.
Note that the corresponding dual fraction $I(S+I)^{-1}$ gives the infected individuals on the boundary and is found in the second equation, in the predation
term. Further, predation on infected prey occurs at rate $w$. The disease is assumed to be unrecoverable, for which the individuals that get it
enter into the class $I$ and can leave it only by being captured by predators, or via natural plus disease-related mortality.
In the last equation the predators' dynamics transpires, which are dependent on the prey for their survival, otherwise they will die at rate $m$.
Predators hunt the healthy and the infected prey alike, but at different rates.

In view of the assumptions stated above, some intrinsic relationships among the parameters hold.
First of all $q \leq w$ and $g \leq f$ since predators hunt infected prey more easily than sound ones; further,
$g<q$ and $f<w$, saying that not the whole captured prey are turned into new predators.

In view of singularities present in (\ref{sistema0}), we need to reformulate the system.

\section{Model reformulation}\label{sec:ref}

We proceed to the singularity elimination, via several steps.
At first, we set $T = \sqrt{S+I}$ in order to remove the square root term.
We thus obtain
\begin{eqnarray}\label{sistema1}
%\begin{cases}
%\displaystyle{
\frac {dS}{dt} &=& \sigma \frac{S^{2}}{T^{2}}-q\frac{PS}{T}+\left( r-\sigma \right)S-\frac{r}{K}ST^{2} 
%}
\\ \nonumber
%\bigskip
%\displaystyle{
\frac {dT}{dt} &=& -\frac{\mu}{2}T-\frac{r}{2K}ST+\left( \frac{r}{2}+\frac{\mu}{2} \right)\frac{S}{T}
+\left( \frac{w}{2}-\frac{q}{2} \right)\frac{PS}{T^{2}}-\frac{w}{2}P
%}
\\  \nonumber
%\displaystyle{
\frac {dP}{dt} &=& -mP+fPT+\left( g-f \right)\frac{PS}{T}.
%}
%\end{cases}
\end{eqnarray}

Then, let $V=ST^{-1}$ in place of $S$.
The system (\ref{sistema1}) becomes
\begin{eqnarray}\label{sistema2}
%\begin{cases}
%\displaystyle{
\frac {dV}{dt} &=& \frac{r}{2K}V^{2}T +\left( \sigma -\frac{r}{2}-\frac{\mu}{2} \right)\frac{V^{2}}{T} +\left( \frac{q}{2}-\frac{w}{2} \right)\frac{PV^{2}}{T^{2}}
\\  \nonumber
&&+\left( r-\sigma+\frac{\mu}{2} \right)V -\frac{r}{K}VT^{2} +\left( \frac{w}{2}-q \right)\frac{PV}{T}
%}
\\  \nonumber
%\bigskip
%\displaystyle{
\frac {dT}{dt} &=& -\frac{\mu}{2}T -\frac{r}{2K}VT^{2} +\left( \frac{r}{2}+\frac{\mu}{2} \right)V +\left(  \frac{w}{2}-\frac{q}{2} \right)\frac{PV}{T} -\frac{w}{2}P %}
\\  \nonumber
%\displaystyle{
\frac {dP}{dt} &=& -mP+fPT+\left( g-f \right)PV.
%}
%\end{cases}
\end{eqnarray}

The third step introduces another new variable, $A=VT^{-1}$ replacing $V$, to reformulate (\ref{sistema2}) as
\begin{eqnarray}\nonumber
%\begin{cases}
%\displaystyle{\dot{A} =
\frac {dA}{dt} &=& \left( \sigma -r -\mu \right) A^{2} +\frac{r}{K}A^{2}T^{2} +\left( q-w \right) \frac{PA^{2}}{T}
+\left( r+ \mu -\sigma \right) A\\ \nonumber
&&-\frac{r}{K}AT^{2} +\left( w-q \right) \frac{PA}{T} \\  \nonumber
%\bigskip
%\displaystyle{\dot{T} =
\frac {dT}{dt} &=& -\frac{r}{2K}AT^{3} -\frac{\mu}{2}T+\left( \frac{r}{2}+\frac{\mu}{2} \right) AT
+\left( \frac{w}{2}-\frac{q}{2} \right)PA -\frac{w}{2}P\\  \label{sistema3}
%\displaystyle{\dot{P} =-mP +ePT}
\frac {dP}{dt} &=& -mP+fPT+\left( g-f \right)PAT.
%\label{sistema3}
%\end{cases}
\end{eqnarray}

This is still unsatisfactory, in view of the presence of the variable $T$ in the denominator.
The next step introduces the variable $U=PT^{-1}$ in place of $P$, to get
the new system with no singularities:
\begin{eqnarray}\label{sistema4}
\frac {dA}{dt} &=& \left( \sigma -r -\mu \right)A^{2} +\frac{r}{K}A^{2}T^{2} +\left( q-w \right) A^{2}U
\\  \nonumber
&&+\left( r+ \mu -\sigma \right)A -\frac{r}{K}AT^{2} +\left( w-q \right)AU, \\ \nonumber
\frac {dT}{dt} &=& -\frac{r}{2K}AT^{3} -\frac{\mu}{2}T+\left( \frac{r}{2}+\frac{\mu}{2} \right)AT
-\frac{w}{2}UT+\left( \frac{w}{2}-\frac{q}{2} \right)AUT, \\ \nonumber
\frac {dU}{dt} &=& \frac{w}{2}U^{2}+ \frac{q-w}{2} AU^{2}+\left(\frac{\mu}{2}-m\right)U+\left( g-f \right)AUT\\ \nonumber
&&- \frac{r+\mu}{2} AU+fUT+\frac{r}{2K}AUT^{2}. 
\end{eqnarray}

Combining all the substitutions made, we find the new variables definitions in terms of the original model variables, as follows
$$
A =\frac{V}{T} =\frac{S}{T^{2}}= \frac{S}{S+I}, \quad U=\frac{P}{T}= \frac{P}{\sqrt{S+I}}, \quad T=\sqrt{S+I},
$$
which allow an interpretation of their meanings.
It follows indeed that
$A$ represents the fraction of healthy prey with respect to the total amount of prey,
$T$ is the total prey population on the edge of the herd and
$U$ denotes the ratio of predators over the total prey population occupying the edge of the area.

\section{Equilibria}\label{sec:equil}

Note first of all that in eliminating singularities we had to divide by $T$, therefore this variable must be different from zero, in fact strictly positive,
so that we exclude possible equilibria with $T=0$. Mathematically, there is a second reason of geometric nature, as $T$ represents the population
of the herd on its boundary, and the latter is certainly never empty for a nonvanishing herd. There are thus only four possible equilibria.

Equilibrium $(A,T,U)=(0,+,0)$ is infeasible since the second equation of (\ref{sistema4}) cannot be satisfied, as it does for
$\left( A,T,U \right) = \left( 0,+,+ \right)$, so that we cannot accept this equilibrium either.

For $\left( A,T,U \right) = \left( +,+,0 \right)$, the first equation of \eqref{sistema4} gives
\begin{equation}\label{3.4}
\frac{r}{K}T^{2}  \left( A-1 \right)= \left( r+ \mu -\sigma \right) \left( A-1 \right)
\end{equation}
so that two cases arise.

If $A = 1$, from the second equation of (\ref{sistema4}) we have $T=\sqrt{K}$, giving the equilibrium
$E_{1}= \left( A_{1}, T_{1}, U_{1}  \right)= \left( 1, \sqrt{K}, 0\right)$ with unconditional feasibility.

Alternatively, if $A<1$, we find
$$
T = \sqrt{\frac{K}{r} \left( r+ \mu -\sigma \right)}, \quad A = \frac{\mu}{\sigma}.
$$
We have thus found the equilibrium
$$
E_{2}=\displaystyle{  \left( A_{2}, T_{2}, U_{2} \right)=
 \left( \frac{\mu}{\sigma}, \sqrt{\frac{K}{r} \left( r+ \mu -\sigma \right)}, 0  \right) }
$$
under the conditions
\begin{equation}\label{E2_feas}
r+\mu -\sigma >0, \quad \mu < \sigma,
\end{equation}
with the second one arising from the very definition of $A$.

{\textbf {Remark.}}
If in $E_{2}$ we let $\mu = \sigma$, we reobtain $E_{1}$.

To find the equilibria with all nonvanishing components $\left( A,T,U \right) = \left( +,+,+ \right)$ that we can call coexistence equilibria,
we sum the second and the third equations of (\ref{sistema4}) to get
\begin{equation}\label{Amm5}
T =\frac m { \left( g-f \right)A + f }, \quad A\ne \frac{f}{f-g}.
\end{equation}
From the first equation of (\ref{sistema4}) we have
$$
\left( A-1 \right)
\left[
\left( \sigma -r -\mu \right) +\frac{r}{K}T^{2}+\left( q-w \right) 
 \right]
=0,
$$
giving again two possibilities.

For $A=1$ we get $T=mg^{-1}$ and the last equation of (\ref{sistema4}) then yields
$$
U=
\frac{r}{g^{2}qK}\left( g^{2}K-m^{2} \right),
$$
which is positive if
\begin{equation}\label{E3_feas}
K > \left( \frac{m}{g} \right)^{2}.
\end{equation}
Thus we found the equilibrium
$$
E_{3}=
\left( A_{3},T_{3},U_{3} \right)
=
\left(  
1, \frac{m}{g},\frac{r}{g^{2}qK}\left( g^{2}K-m^{2} \right)
\right)
$$
with feasibility condition (\ref{E3_feas}).

If instead $A\ne 1$ we solve the system
\begin{eqnarray}\label{4.3}
\left( \sigma -r -\mu \right) +\frac{r}{K}T^{2}+\left( q-w \right)U=0,\\ \nonumber
-\frac{r}{2K}AT^{2} -\frac{\mu}{2}+\left( \frac{r+\mu}{2} \right)A -\frac{w}{2}U+\left( \frac{w-q}{2} \right)AU=0,
%\\
%T=\frac{m}{\left( g-f \right)A + f }.
\end{eqnarray}
with $T$ given by the first equation in (\ref{Amm5}).
Now in the first equation \eqref{4.3} write $U$ as a function of $A$:
\begin{equation}\label{4.3.a}
\left( q-w \right)U = \left( r+\mu -\sigma \right) -\frac{r}{K} \frac{m^{2}}{\left[ \left(g-f\right)A+f \right]^{2} }.
\end{equation}

Now if $q-w=0$
the first equation of (\ref{4.3})
simplifies to give
\begin{equation}\label{T5}
T= \sqrt{\frac{K}{r}\left( r +\mu -\sigma \right) },
\end{equation}
provided $r+\mu -\sigma > 0$, an assumption that we are making from now on.
Substituting into the second equation (\ref{4.3}) we find
$$
U= \frac{\sigma A -\mu}{w},
$$
which is nonnegative if $ A\ge \mu \sigma^{-1}$.
From the first equation in (\ref{Amm5}) we then obtain
$$
A =\frac{1}{ g-f } \left[ m\sqrt{\frac{r}{K\left( r +\mu -\sigma \right)}}-f\right] .
$$
Recalling the assumption $g<f$, $A$ will be nonnegative if and only if
$$
K>\left( \frac{m}{f} \right)^{2} \frac{r}{ r +\mu -\sigma} .
$$
We finally have the explicit expression of $U$ as follows,
$$
U=\frac{\sigma}{w} \frac{1}{g-f} 
\left[
m\sqrt{\frac{r}{K\left( r +\mu -\sigma \right)}}-f
\right]
-\frac{\mu}{w},
$$
which is nonnegative when $A> \mu\sigma^{-1}$, i.e. for
\begin{equation}\label{4.5}
m\sqrt{\frac{r}{K\left( r +\mu -\sigma \right)}}<\frac{\mu}{\sigma}\left( g-f \right)+f.
\end{equation}
The right hand side is positive if $\left(\sigma -\mu\right)f+\mu g>0 $ and from this the above restriction can be rewritten as
$$
K>\left[
\frac{m \sigma}{ \left(\sigma -\mu\right)f+\mu g }
\right]^{2}
\frac{r}{r+\mu -\sigma}.
$$

In summary we found the equilibrium
$E_{4}=\left( A_{4},T_{4},U_{4} \right)$
where, explicitly,
\begin{eqnarray*}
A_{4}=\frac{1}{ g-f } \left[ m\sqrt{\frac{r}{K\left( r +\mu -\sigma \right)}}-f\right], \quad
T_{4}=\sqrt{\frac{K}{r}\left( r +\mu -\sigma \right) },\\
U_{4}=\frac{\sigma}{ w\left( g-f\right)  } \left[m\sqrt{\frac{r}{K\left( r +\mu -\sigma \right)}}-f\right]-\frac{\mu}{w},
\end{eqnarray*}
with feasibility conditions $r+\mu-\sigma >0$, $q=w$ and
\begin{eqnarray*}
K>\left[ \frac{m \sigma}{ \left(\sigma -\mu\right)f+\mu g}
\right]^{2} \frac{r}{r+\mu -\sigma}, \quad
\sigma f+\mu g>\mu f, \quad
K> \frac{m^2r}{ f^2(r +\mu -\sigma)}. \\
\end{eqnarray*}
These conditions can be simplified, observing that from $g-f<0$ it follows
$$
\frac{m\sigma}{f\sigma}<\frac{m\sigma}{\left(g-f\right)\mu+f\sigma }
$$
so that if $\left(\sigma -\mu\right)f+\mu g>0$,
$$
K>\left( \frac{m}{f} \right)^{2} \frac{r}{ r +\mu -\sigma}
$$
is implied by the condition
$$
K>\left( \frac{m \sigma}{\left[ \left(\sigma -\mu\right)f+\mu g \right]} \right)^{2}
\frac{r}{r+\mu -\sigma}.
$$
Thus feasibility conditions for $E_{4}$ become just the following ones
\begin{eqnarray}\label{E4_feas}
K> \frac{m^2 r \sigma^2}{\left[ \left(\sigma -\mu\right)f+\mu g \right]^2 (r+\mu -\sigma)} ,
\quad
\sigma f+\mu g>\mu f,
\quad
r+\mu>\sigma ,
\quad
q=w.
\end{eqnarray}

We now address the case $q-w<0$.
In this situation from \eqref{4.3.a} we find
\begin{equation}\label{U5}
U=
\frac{1}{\left(q-w \right) }
\left[
\left( r+\mu -\sigma \right) -\frac{r}{K} \frac{m^{2}}{ \left[ \left(g-f\right)A+f \right]^{2}  }
\right]
\end{equation}
and substituting the values of $T$ and $U$ into the second equation of \eqref{4.3}, we obtain
\begin{equation}\label{4.6}
-\mu - \frac{w\left( r+\mu -\sigma \right)}{q-w }
+\frac{rwm^{2}}{ \left[ K\left(g-f\right)A+f \right]^{2}  \left( q-w \right)}
+\sigma A=0.
\end{equation}
From this, with some algebra, we are led to the following cubic equation for $A$:
\begin{eqnarray}\label{cubic_E5}
P\left(A\right) \equiv \sum_{k=0}^3 b_k A^k =0,
\end{eqnarray}
where
% $b_3=1$, CHIARIRE SE QUESTO E' VERO E CONTROLLARE INDICI QUI SOTTO!!!!!!
\begin{eqnarray*}
b_3= \sigma \left( q-w \right) K^2\left(g-f\right)^{2},\\
b_2 = 
2f\left(g-f\right)\sigma \left( q-w \right)K
%\right]A^{2}
%+\left[
+\left( w\sigma -rw-q\mu \right)K^2\left(g-f\right)^{2},\\
b_1= f^{2}\sigma \left( q-w \right)K +2f\left(g-f\right)\left( w\sigma -rw-q\mu \right)K,\\
b_0=f^{2}\left( w\sigma -rw-q\mu \right)K +m^{2}rw.
\end{eqnarray*}

%\newpage
It has always a real root, and we now seek sufficient conditions for a nonnegative real root. Since $q < w$, it follows that
\begin{equation}\label{infinito}
\lim_{A \to +\infty} P\left(A\right) = -\infty,
\end{equation}
so that if the constant term is positive, at least one positive real root must exist. This occurs if
\begin{equation}\label{4.8}
f^{2}\left( w\sigma -rw-q\mu \right)K + m^{2}rw>0,
\end{equation}
which is trivial in case 
\begin{equation}\label{E5_feas1}
w\sigma -rw-q\mu \geq 0,
\end{equation}
otherwise it leads to
\begin{equation}\label{E5_feas2}
w\sigma -rw-q\mu < 0, \quad
K< \left( \dfrac{m}{f} \right)^{2} \dfrac{rw}{rw+q\mu -w\sigma}.
\end{equation}
In summary the equilibrium $E_{5}=\left(A_{5},T_{5},U_{5} \right)$ arises with first component given by the positive root of (\ref{cubic_E5})
and the remaining ones by (\ref{T5}) and (\ref{U5}), which need to be nonnegative, and further feasibility conditions given by (\ref{E5_feas1}) or
(\ref{E5_feas2}).

\section{Stability}\label{sec:stab}
The elements of the Jacobian matrix $J=(J_{ik})$, $i,k=1,2,3$ are
\begin{eqnarray*}
J_{11}= (2A-1)[\left( \sigma -r -\mu \right) +\left(q-w \right)U +\frac{r}{K}T^{2} ] \quad
%J_{11}= (2A-1)[\left( \sigma -r -\mu \right) +\left(q-w \right)U] +\frac{r}{K}T^{2} (2A-1) \\
%%-\frac{r}{K}T^{2} \\
%J_{11}= (2A-1)[\left( \sigma -r -\mu \right) +\left(q-w \right)U] +\frac{2r}{K}AT^{2} 
%-\frac{r}{K}T^{2} \\
%J_{11}= \left( \sigma -r -\mu \right)(2A-1) +\frac{2r}{K}AT^{2} +\left(q-w \right)U (2A-1) 
%-\frac{r}{K}T^{2} \\
%J_{11}= 2\left( \sigma -r -\mu \right)A +\frac{2r}{K}AT^{2} +2\left(q-w \right)AU + r +\mu -\sigma
%-\frac{r}{K}T^{2} +\left( w-q\right)U\\
J_{12}= 2\frac{r}{K}AT\left( A-1 \right)\\
J_{13}= \left( q-w \right) A \left( A-1\right)\quad
J_{21}=-\frac{r}{2K}T^{3}+ \frac{r+\mu}{2} T + \frac{w-q}{2} UT\\
J_{22}= -\frac{3r}{2K}AT^{2}-\frac{\mu}{2}+ \frac{r+\mu}{2} A-\frac{w}{2}U + \frac{w-q}{2} AU\quad
J_{23}=-\frac{w}{2}T+ \frac{w-q}{2} AT\\
J_{31}= \frac{q-w}{2} U^{2} +\left( g-f \right)UT - \frac{r+\mu}{2} U +\frac{r}{2K}UT^{2}\\
J_{32}= \left( g-f\right)AU +fU +\frac{r}{K}AUT\\
J_{33}=wU+\left( q-w\right)AU + \frac{\mu}{2}-m +\left( g-f \right)AT - \frac{r +\mu}{2} A +fT+\frac{r}{2K}AT^{2}
\end{eqnarray*}
Observe that since $A =S(S+I)^{-1}\leq 1$ and $q<w$ two of the above terms have a fixed sign:
$$
J_{12} \leq 0, \quad J_{13} \geq 0.
$$

The Jacobian's eigenvalues at $E_1$ are
$\lambda _{1} =\sigma - \mu$, $\lambda _{2} =-r$, $\lambda _{3} =-m +g\sqrt{K}$, from which the stability conditions follow
\begin{equation}\label{E1_stab}
\dfrac{\mu}{\sigma} >1,\quad K < \frac{m^2}{g^2} .
\end{equation}

The Jacobian at $E_2$ gives one eigenvalue as
$$
\lambda _{1} =
\sqrt{\frac{K}{r} \left( r+ \mu -\sigma \right)}\left[f+\left(g-f\right)\frac{\mu}{\sigma} \right]-m.
$$
from which the stability condition follows
\begin{equation}\label{E2_stab}
K<\left[\frac{m\sigma}{\left( \sigma - \mu \right)f+g\mu }\right]^{2}
\frac{r}{r+ \mu -\sigma }
\end{equation}
having used the fact that $\left( \sigma - \mu \right)f+g\mu>0$ and the first condition (\ref{E2_feas}).
% $r+ \mu -\sigma >0$.
The other two eigenvalues are the roots of\\[2ex]
\begin{equation}\label{S1g}
\lambda ^{2}+\frac{\mu}{\sigma} \left( r+ \mu -\sigma \right)\lambda +\mu \left( 1-\frac{\mu}{\sigma} \right)\left( r+ \mu -\sigma \right)=0.
\end{equation}
In view of the feasibility conditions (\ref{E2_feas}), the Routh-Hurwitz stability conditions for (\ref{S1g}) hold.
Stability of $E_2$ is therefore regulated only by (\ref{E2_stab}).

At $E_3$ again one eigenvalue is immediate,
$$
\lambda _{1}= \displaystyle{\left( \sigma-\mu \right) +\frac{rw}{g^{2}qK}\left( m^{2}-g^{2}K \right)}.
$$
It is negative if and only if 
$g^{2}K\left[ q\left( \sigma -\mu \right) -rw \right] < -m^{2}rw $.
But this cannot happen if
$ q\left( \sigma -\mu \right) -rw \geq 0$. Conversely, we are lead to the stability conditions
\begin{equation}\label{E3_stab1}
K > \frac{rwm^2}{g^2(rw+q\mu -q\sigma)} , \quad  rw + q \sigma > q\mu .
\end{equation}
The other eigenvalues come from the quadratic
\begin{equation}\label{E3_char}
\lambda ^{2} 
+\frac{r}{2g^{2}K} \left( 3m^{2}-g^{2}K \right)\lambda
+\frac{mr}{2g^{4}K} \left( g^{2}K-m^{2} \right)
\left(2mr+g^{2}K\right) =0.
\end{equation}
From the (strict) feasibility conditions (\ref{E3_feas}) for $E_{3}$, the constant term is always positive. Imposing that also the coefficient
of the linear term is positive, we obtain the second stability condition,
\begin{equation}\label{5.2}
K <3 \frac{m^2}{g^2} .
\end{equation}
In summary, $E_3$ is feasible and stable for
\begin{equation}\label{E3_stab}
0< \max\left\{ 1, \frac {rw}{rw+q\mu -q\sigma}\right\} < K \frac {g^2}{m^2} <3.
\end{equation}

For the equilibrium $E_4$ some of the Jacobian entries, in view of the feasibility conditions (\ref{E4_feas}) have fixed signs, as follows
\begin{eqnarray*}
J_{4_{12}}=
\frac{2r}{K}\sqrt{\frac{K}{r}\left( r +\mu -\sigma \right)}\frac{1}{ g-f } \left[ m\sqrt{\frac{r}{K\left( r +\mu -\sigma \right)}}-f\right]\\
\cdot \left[
\frac{1}{ g-f } \left( m\sqrt{\frac{r}{K\left( r +\mu -\sigma \right)}}-f\right)-1
\right]<0,\\
J_{4_{21}}=
\frac{\sigma}{2}\sqrt{\frac{K}{r}\left( r +\mu -\sigma \right)}>0, \quad
J_{4_{23}}=
-\frac{w}{2}\sqrt{\frac{K}{r}\left( r +\mu -\sigma \right)}<0,\\
J_{4_{22}}=
-\frac{ r+\mu -\sigma }{g -f }\left[ m\sqrt{\frac{r}{K\left( r +\mu -\sigma \right)}}-f\right]<0,\\
J_{4_{32}} 
=
\left\lbrace
\frac{\sigma}{ w\left( g-f\right)  } \left[m\sqrt{\frac{r}{K\left( r +\mu -\sigma \right)}}-f\right]-\frac{\mu}{w}
\right\rbrace \\
\times \left[
m\sqrt{\frac{r}{K\left( r +\mu -\sigma \right)}}
%\right.\\
+
%\left.
%\frac{1}{ \left( g-f \right)} \left( m\sqrt{\frac{r}{K\left( r +\mu -\sigma \right)}}-f\right)
%\sqrt{\frac{r}{K}\left( r +\mu -\sigma \right)}
%\right]>0,\\
\frac{1}{ g-f} \left(m\frac rK-f
\sqrt{\frac{r}{K}\left( r +\mu -\sigma \right)}\right)
\right]>0,
\end{eqnarray*}
while the remaining two must agree, since the same factor appears in the two elements, although the sign is not decided:
\begin{eqnarray*}
J_{4_{31}}
=
\frac{1}{w}
\left[
\sigma\left( m\sqrt{\frac{r}{K\left( r +\mu -\sigma \right)}}-f \right)-\mu\left(g-f\right) 
\right]
\cdot
\\
\left[
\sqrt{\frac{K}{r}\left( r +\mu -\sigma \right)} -\frac{\sigma}{2\left(g-f \right) }
\right],
\\
J_{4_{33}}
=\frac{1}{2\left( g-f \right)}
\left[
\sigma\left( m\sqrt{\frac{r}{K\left( r +\mu -\sigma \right)}}-f \right)-\mu\left(g-f\right) 
\right].
\end{eqnarray*}
We now study the signs of $J_{4_{31}}$ and $J_{4_{33}}$.
Considering $J_{4_{31}}$ and using the feasibility condition (\ref{E4_feas}) we find that for its positivity we must have
$$
\frac{r}{K\left( r +\mu -\sigma \right)}
> \left[\frac{\mu\left(g-f\right) +f\sigma}{m\sigma}\right]^{2}.
$$
But this contradicts the feasibility condition (\ref{E4_feas}), so that it must be negative. In summary we then have
$$
J_{4_{31}}<0, \quad J_{4_{33}}<0.
$$
Thus the resulting structure of the Jacobian matrix is
$$
J_{4} =
\left(
\begin{array}{ccc}
 \bigskip 
0 & - & 0\\
 \bigskip 
+ & - & - \\
- & + & -
\end{array}
\right)
\equiv
\left(
\begin{array}{ccc}
\bigskip
0 & Z & 0\\
\bigskip
B & C & D \\
E & F & G
\end{array}
\right).
$$
The characteristic equation is now a cubic,
\begin{equation}\label{7.1}
\sum_{k=0}^3 a_k \lambda^k\equiv 
\lambda^{3}
-\left( C+G \right)\lambda ^{2}
- \left( ZB+FD-CG \right) \lambda
-Z \left( ED-BG \right)=0.
\end{equation}
Using the signs of $Z$, $B$, $C$, $D$, $E$, $F$ and $G$ all the coefficients $a_{k}$, $k=0,\ldots,3$ are positive.
We can thus use the Li\'enard-Chipart criterion, a particular case of the Routh-Hurwitz criterion,
thereby determining the sign of the eigenvalues imposing that the following determinant be positive:
\begin{eqnarray}\nonumber
D_{2} =
\left|
\begin{array}{cc}
\bigskip
a_{2} & a_{0}\\
a_{3} & a_{1}
\end{array}
\right|
=
\left|
\begin{array}{cc}
\bigskip
-\left( C+G \right) & -Z \left( ED-BG \right)\\
1 & - \left( ZB+FD-CG \right)
\end{array}
\right| =\\ \label{E4_stab}
=\left( C+G \right) \left( ZB+FD-CG \right) +Z \left( ED-BG \right)>0.
\end{eqnarray}
%For stability we need that $D_{2}>0$, i.e.,
%\begin{displaymath}
%\left( C+G \right) \left( ZB+FD-CG \right) +Z \left( ED-BG \right)>0.
%\end{displaymath}
We can conclude for this case that $E_{4}$ is stable if (\ref{E4_stab}) holds.
%\begin{eqnarray}\label{E4_stab}
%K>\left( \frac{m \sigma}{\left[ \left(\sigma -\mu\right)f+\mu g \right]} 
%\right)^{2}
%\frac{r}{r+\mu -\sigma}, \quad
%\left(\sigma -\mu\right)f+\mu g>0
%\quad 
%r+\mu-\sigma >0
%\quad
%q=w \quad
%\left( C+G \right) \left( ZB+FD-CG \right) +Z \left( ED-BG \right)>0
%\end{eqnarray}

Stability of $E_5$ is investigated numerically.

\section{Bifurcations}\label{sec:Hopf}
Note that transcritical bifurcations further arise between $E_{1}$ and $E_{2}$, compare (\ref{E1_stab}) with (\ref{E2_feas})
and the remark below it, as well as $E_{1}$ and $E_{3}$,
see (\ref{E1_stab}) and (\ref{E3_feas}).

We then try to establish if there are special parameter combinations for which Hopf bifurcations arise.
For this purpose, we need purely imaginary eigenvalues. This is easy to assess for a quadratic characteristic equation,
$\lambda^{2} + b \lambda +c =0$
since we need the linear term to vanish, $b=0$, and the constant term to be negative, $c<0$. For a generic cubic
of the form
\begin{equation}\label{E1g}
a_{3}\lambda ^{3}+a_2\lambda ^{2} +a_{1}\lambda +a_{0}=0 
\end{equation}
instead, we need the following condition
\begin{displaymath}
a_{1}a_{2}-a_{0}=0.
\end{displaymath}
Clearly at $E_1$ no bifurcation arises, since the eigenvalues are all real.
At $E_2$ we need
\begin{equation}\label{E2}
b=\displaystyle{\frac{\mu}{\sigma} \left( r+ \mu -\sigma \right)=0},\quad
c=\displaystyle{\mu \left(1- \frac{\mu}{\sigma} \right)\left( r+ \mu -\sigma \right)}>0
\end{equation}
but these conditions contradict each other. We conclude that at $E_{2}$ no Hopf bifurcations can carise.

At $E_{3}$ the characteristic equation factors, and the quadratic (\ref{E3_char}) 
from feasibility (\ref{E3_feas}) has a positive constant term. Imposing that the linear term vanishes, we find the value
\begin{equation}\label{E3_Hopf}
K^{\dagger} \equiv 3 \frac{m^2}{g^2}
\end{equation}
for which a limit cycle appears.
%We will show this situation and investigate the remaining ones for $E_4$ and $E_5$ with simulations.
%In the study of the stability, focusing on the bifurcations, we discovered that, around $E_{3}$, limit cycles should arise
%when $K$ crosses the threshold value $K^{\dagger}=3 m^2g^{-2}$.
In Figure \ref{f:5} we present
a simulation of the two-dimensional limit cycle for the parameter values
$\sigma=0.5$, $r=0.5,\mu=0.4$, $q=0.2$, $w=0.5$, $m=0.2$, $g=0.1$, $f=0.3$.
Oscillations appear only for the second and the third variables, while for first one remains
at the fixed level $A=1$, to mean that the system is disease-free.
Predators survive together with the healthy individuals but with persistent oscillations of the two populations.
%In Figure \ref{f:6} a three-dimensional phase-space portrait of the limit cycle is given.
Note that this bifurcation is due to demographics effects only, not to epidemiological ones, compare with the earlier works \cite{AV,APV},
since it occurs for the same parameter value, with changes in notations only.

%It follows that the number of prey on the edge is given by the number of healthy prey on it. The third section of the figure tells us that predators survive together with the healthy individuals.

%\newpage

\begin{figure}[h]
\centering
\includegraphics[scale=0.47]{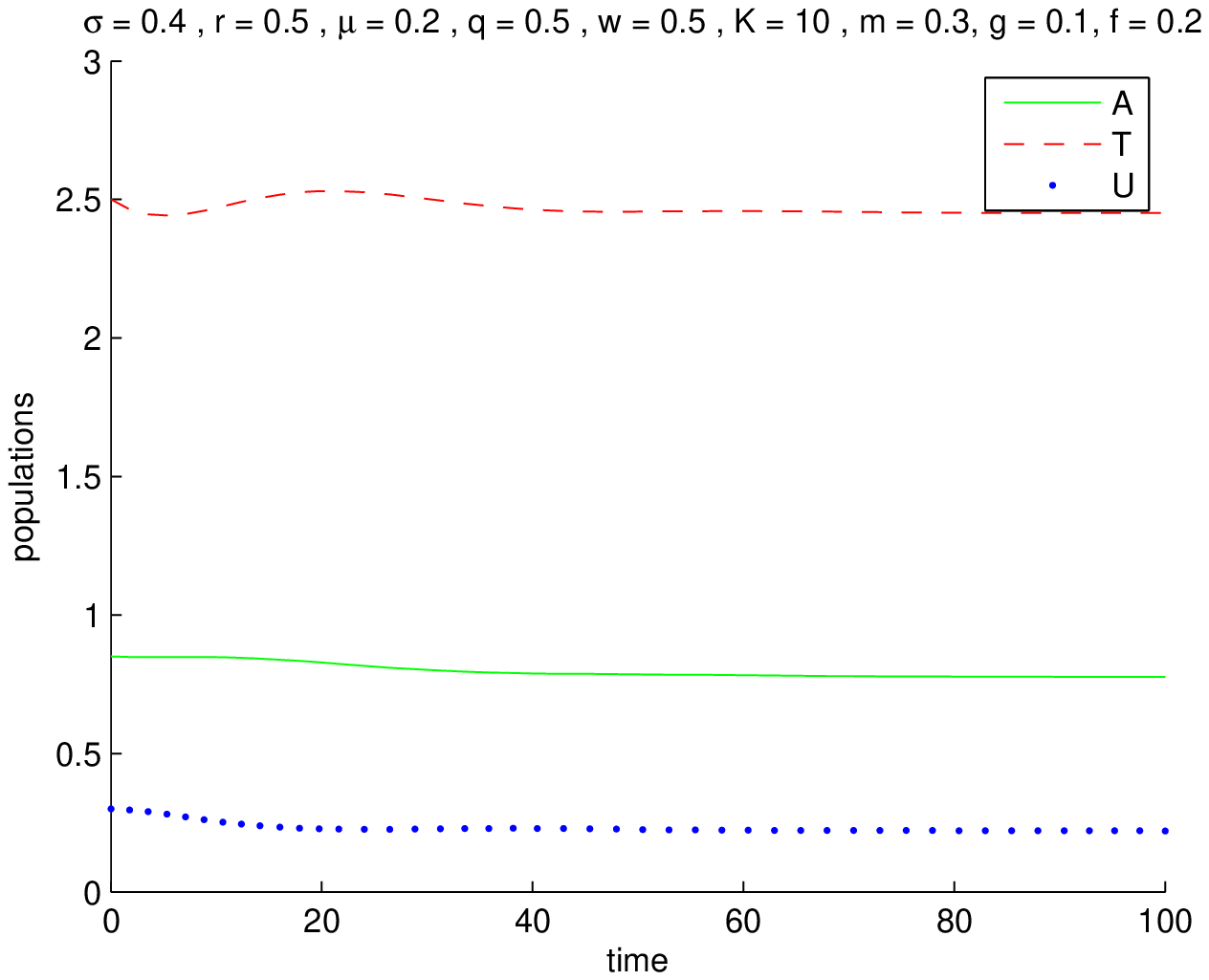}
\includegraphics[scale=0.47]{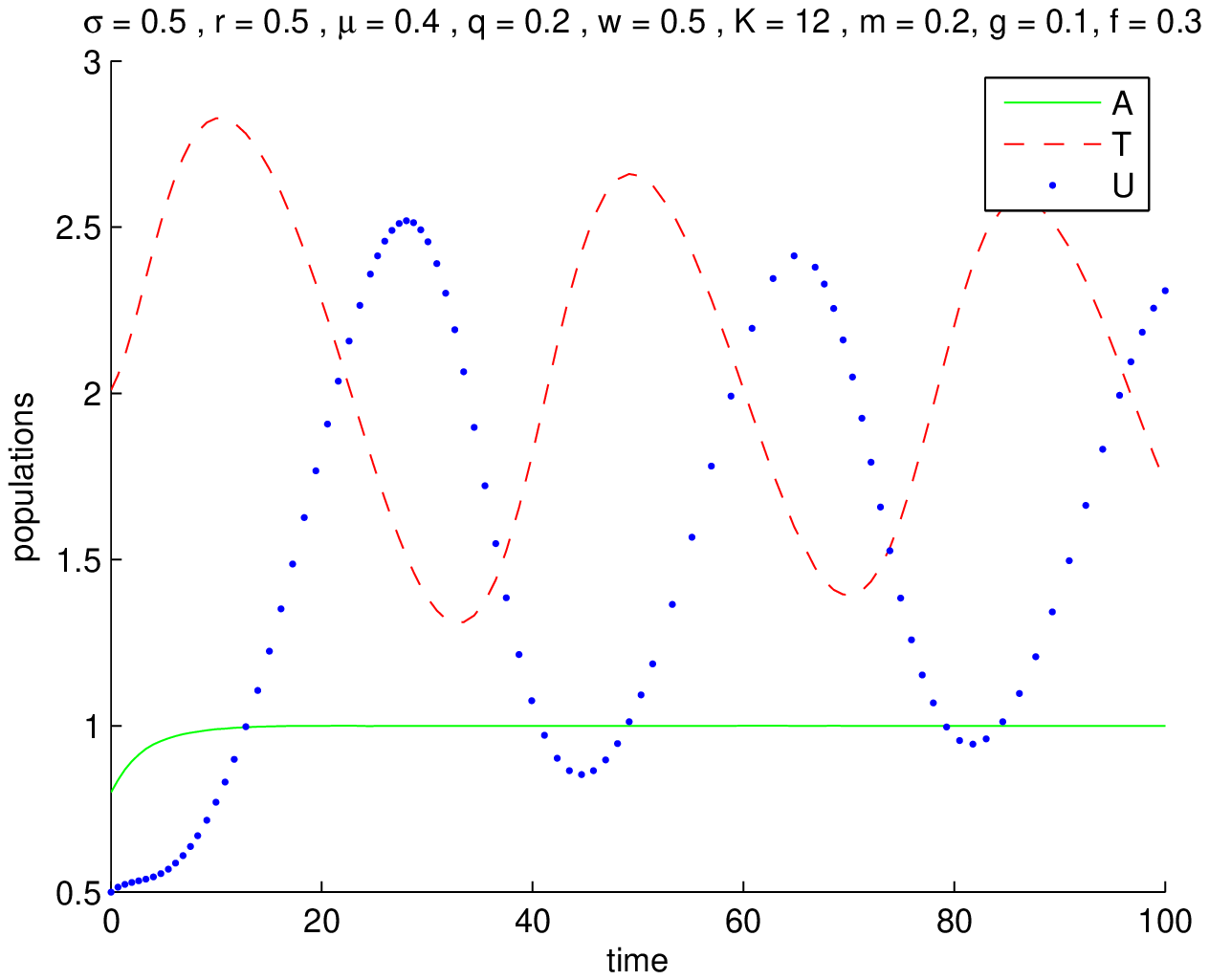}
\caption{Left: Stable equilibrium $E_4$ for the choice $\sigma=0.4$, $r=0.5$, $\mu=0.2$, $q=0.5$, $w=0.5$, $K=10$, $m=0.3$, $g=0.1$, $f=0.2$;
Two-dimensional limit cycle around $E_4$ obtained for the parameter values
$\sigma=0.5$, $r=0.5,\mu=0.4$, $q=0.2$, $w=0.5$, $m=0.2$, $g=0.1$, $f=0.3$, for which $K=K^{\dagger}\equiv 3\left(\dfrac{m}{g}\right)^2.$}
\label{f:5}
\end{figure}

\section{Interpretation of the system's evolution}\label{sec:interp}

At equilibrium $E_{1}$, we have $U_{1}=0$ so that $P_1=0$ and the predators vanish. Further, $A_1=1$ implying that $I_1=0$. Thus only healthy prey survive,
at the environment's carrying capacity, due to the model assumption of logistic growth, $T_1=\sqrt K$ indeed implies in this case $S_1=K$.
The equilibrium $E_{1}$ thus represents the situation where the only population which survives in the habitat is represented by the healthy prey.
The fact that infected individuals are extinguished is consistent with the stability conditions of $E_{1}$. In fact,
the latter require that the disease incidence be lower than the disease-related mortality rate. Thus infected individuals die
faster than they are recruited and ultimately there are not enough infectious individuals to propagate the disease.
Its stable behavior can be obtained
%is shown in Figure \ref{f:1} 
for the parameter values
$\sigma =0.2$, $r=0.5$, $K=5$, $\mu =0.4$, $q=0.2$, $w=0.5$, $m=0.8$, $g=0.1$, $f=0.3$.

At $E_2$ the request that $A<1$ tells us that neither healthy nor infected prey disappear from the system,
while, as in the previous case, all the predators die since $U_2=0$. Therefore
the disease remains endemic among the prey, while predators do not survive. Note once again that the the point $E_{2}$ becomes equilibrium $E_{1}$
if we assume that the disease transmission rate equals the disease mortality rate. In such case the disease can thus be eradicated.
At the equilibrium $E_2$ predators get extinguished, but the disease remains endemic. In this situation the opposite condition
of equilibrium $E_1$ must be verified, namely the disease-related mortality rate is lower than the disease incidence.
This suggests that it is reasonable to expect that the population of infected prey survives.
%Figure \ref{f:2} contains a simulation 
A set of parameter values leading to this equilibrium is for instance given by
$\sigma=0.5$, $ r=0.5$, $K=5$, $\mu=0.4$, $q=0.2$, $w=0.5$, $m=0.8$, $g=0.1$, $f=0.3.$.

At $E_{3}$
we have again that $A_{3}=1$, so that
$I=0$
and in this case the disease gets eradicated from the ecosystem, while the predators and healthy prey survive together.
This is the only equilibrium for which we have proved analytically the existence of bifurcations,
for the particular value of the prey carrying capacity $K^{\dagger}=3 m^2g^{-2}$.
This third equilibrium $E_3$ is stably achieved e.g. for the parameter set
%shows coexistence of predators and healthy prey, with the disease eradicated.
%Figure \ref{f:3} shows it for the parameter values
$\sigma =0.5$, $r=0.5$, $K=5$, $\mu= .4$, $q=0.2$, $w=0.5$, $m=0.2$, $g=0.1$, $f=0.3.$.

At $E_4$ and $E_5$ we have coexistence, with the point $E_3$ being a particular case of the latter equilibria, when $A=1$. Further
$E_4$ and $E_5$ differ because in the first case $q=w$, i.e. the infected and healthy prey are hunted at the same rate by predators,
and therefore it can be regarded as a special case of $E_5$. As for the latter, note that for the particular situation in which
$f^{2}K \left( rw +q\mu -w\sigma \right)=m^{2}rw$ we find $A_5$, as the cubic (\ref{cubic_E5}) goes through the origin.
This implies that the healthy prey are wiped out. Therefore in this situation the ecosystem thrives, with predators and only infected prey.
The equilibrium $E_{4}$ can be obtained by the choice
%represents the possibility that all the populations in the system survive, i.e. $E_{4}$ represents the coexistence equilibrium.
%The simulation of Figure \ref{f:4} shows it for the following parameter values:
$\sigma =0.4$, $r=0.5$, $\mu =0.2$, $q=0.5$, $w=0.5$, $m=0.3$, $f=0.2$, $g=0.1$, $K=10.$.
Instead, for the equilibrium $E_5$ our extensive simulations seem to indicate its instability.

\section{The poisonous prey}\label{pois}

Here we consider the situation in which ingested infectious prey harm predators.
\begin{eqnarray}\label{sistema0poisonous}
%\begin{cases}
%\displaystyle{
\frac {dS}{dt} &=& rS\left( 1-\frac{S+I}{K} \right)-\sigma \frac{SI}{S+I}-qPS \frac{\sqrt{S+I}}{S+I} %}
\\ \nonumber
%\bigskip
%\displaystyle{
\frac {dI}{dt} &=& \sigma \frac{SI}{S+I}-w PI \frac{\sqrt{S+I}}{S+I}- \mu I %}
\\ \nonumber
%\displaystyle{
\frac {dP}{dt} &=& -mP+gP\frac{S}{\sqrt{S+I}} -fP\frac{I}{\sqrt{S+I}} 
%-mP +eP\sqrt{S+I}, %}
%\end{cases}
\end{eqnarray}

If it is more difficult to capture healthy animals, we need the following
assumptions on the parameters $g\le q$, $f\le w$, $q\le w$.

Introducing the same new parameters as for model (\ref{sistema0}), we are then led to the singularity-free system
\begin{eqnarray}\label{sistema5}
\frac {dA}{dt} &=& \left( \sigma -r -\mu \right)A^{2} +\frac{r}{K}A^{2}T^{2} +\left( q-w \right) A^{2}U
\\  \nonumber
&&+\left( r+ \mu -\sigma \right)A -\frac{r}{K}AT^{2} +\left( w-q \right)AU, \\ \nonumber
\frac {dT}{dt} &=& -\frac{r}{2K}AT^{3} -\frac{\mu}{2}T+\left( \frac{r}{2}+\frac{\mu}{2} \right)AT
-\frac{w}{2}UT+\left( \frac{w}{2}-\frac{q}{2} \right)AUT, \\ \nonumber
\frac {dU}{dt} &=& \frac{w}{2}U^{2}+ \frac{q-w}{2} AU^{2}+\left(\frac{\mu}{2}-m\right)U+\left( g+f \right)AUT\\ \nonumber
&&- \frac{r+\mu}{2} AU-fUT+\frac{r}{2K}AUT^{2}. 
\end{eqnarray}

The equilibria, here denoted by $P_i=(A_i,T_i,U_i)$, again contain $T\neq 0$, and are as follows.
$P_1=(1,\sqrt K,0)$,
$$
P_2=\left( \frac {\mu}{\sigma},\sqrt {\frac Kr(r+\mu-\sigma)},0\right), \quad
P_3=\left( 1,\frac mg, \frac r{g^2qK}(g^2K-m^2)\right)
$$ 
Since $A\le 1$, and in view of the root in its expression, $P_2$ requires for feasibility 
\begin{equation}\label{feas_P2}
r+\mu\ge \sigma \ge \mu.
\end{equation}
Clearly $P_1$ is the limiting case of $P_2$ when $\mu=\sigma$.
$P_3$ is feasible for 
\begin{equation}\label{feas_P3}
K\ge \frac {m^2}{g^2}.
\end{equation}
Assuming $q=w$, with calculations that mimic those of the first part, also the equilibrium $P_4$ can be found, with
\begin{equation}\label{P4_pois}
A_4=\frac {m+f T_4}{(f+g)T_4}, \quad
T_4=\sqrt {\frac Kr(r+\mu-\sigma)}, \quad
U_4=\frac {\sigma(m+fT_4)}{w(f+g)T_4}-\frac {\mu}w.
\end{equation}
The feasibility conditions for $P_4$ are
\begin{equation}\label{feas_P4}
r+\mu\ge \sigma , \quad q=w, \quad T_4g\le \frac {\sigma}{\mu} \left[ m+fT_4 \right]-fT_4.
\end{equation}
Finally, for $q<w$, we can establish the existence of the equilibrium $P_5$ proceeding as done for the formulae
(\ref{U5}), with $T_{5}\equiv T_4$ and obtaining now
\begin{equation}\label{U5_pois}
U_5=
\frac{1}{q-w}
\left[
r+\mu -\sigma -\frac{r}{K} \frac{m^{2}}{ \left[ \left(g+f\right)A_5-f \right]^{2}  }
\right]
\end{equation}
and (\ref{4.6}), to find for $A_5$ again a cubic like (\ref{cubic_E5}), but with coefficients this time given by
\begin{eqnarray*}
b_3= \sigma \left( q-w \right) K\left(g+f\right)^{2},\\
b_2 = 
\left( w\sigma -rw-q\mu \right)K^2\left(g+f\right)^{2}
-2f\left(g+f\right)\sigma \left( q-w \right)K,\\
b_1= f^{2}\sigma \left( q-w \right)K -2f\left(g+f\right)\left( w\sigma -rw-q\mu \right)K,\\
b_0=f^{2}\left( w\sigma -rw-q\mu \right)K +m^{2}rw.
\end{eqnarray*}
In view of the assumption $q< w$, again we find the behavior expressed by (\ref{infinito}), so that it is enough to require $b_0>0$
to have a positive root. It is given once again by (\ref{4.8}), 
which is trivial if (\ref{E5_feas1}) holds, 
otherwise again it leads to (\ref{E5_feas2}).
$P_5$ is then feasible if, recalling (\ref{P4_pois}), we have $T_{5}\equiv T_4 \ge 0$, which amounts to the first one of
(\ref{feas_P4}), and if also $U_{5}\ge 0$. From (\ref{U5_pois}), this imposes a lower bound on $A_5$. In summary feasibility is obtained for
\begin{equation}\label{feas_P5}
q<w, \quad r+\mu\ge \sigma, \quad A_{5}\ge \frac 1{f+g}\left[\frac {m\sqrt r}{\sqrt{K(r+\mu-\sigma)} }+f \right].
\end{equation}

The Jacobian of system (\ref{sistema5}) differs from the one relative to model (\ref{sistema0}) only in some signs of the elements of
the last row, namely
\begin{eqnarray*}
J_{31}= \frac{q-w}{2} U^{2} +\left( g+f \right)UT - \frac{r+\mu}{2} U +\frac{r}{2K}UT^{2}\\
J_{32}= \left( g+f\right)AU -fU +\frac{r}{K}AUT\\
J_{33}=wU+\left( q-w\right)AU + \frac{\mu}{2}-m +\left( g+f \right)AT - \frac{r +\mu}{2} A -fT+\frac{r}{2K}AT^{2}
\end{eqnarray*}

For equilibrium $P_1$ the eigenvalues are easily explicitly obtained, to give the stability conditions
\begin{equation}\label{stab_P1}
\sigma < \mu, \quad K<  \frac {m^2}{g^2} .
\end{equation}
It can be numerically shown to be stable with the same parameter choice employed for corresponding equilibrium $E_1$ in the harmless situation.

For $P_2$, one eigenvalue is explicit, and then we obtain a quadratic characteristic equation for which the Routh-Hurwitz conditions hold
unconditionally, in view of the feasibility condition (\ref{feas_P2}). Stability is then achieved if
$$
\sqrt {\frac Kr(r+\mu-\sigma)}\left[ \frac {\mu}{\sigma} (f+g) -f\right] <m
$$
which gives one of the two alternative following sets of conditions
\begin{eqnarray}\label{stab_P2a}
g\le \frac f{\mu} (\sigma - \mu); \\ \label{stab_P2b}
g> \frac f{\mu} (\sigma - \mu), \quad K<  \frac r{r+\mu-\sigma} \frac {m^2\sigma^2}{\left[ g\mu-f(\sigma-\mu) \right]^2}.
\end{eqnarray}
This equilibrium arises for the same parameter values used for $E_2$.

At $P_3$, again the Jacobian factors to give one explicit eigenvalue and a quadratic. The Routh-Hurwitz condition on the latter
reduces only to requiring
\begin{eqnarray}\label{stab_P3a}
K< 3 \frac {m^2} {g^2}
\end{eqnarray}
while for the negativity of the former we need
\begin{eqnarray}\label{stab_P3b}
rw+q\mu-q\sigma>0, \quad K> \frac {m^2rw}{g^2(rw+q\mu-q\sigma)}.
\end{eqnarray}
This equilibrium can be obtained with the same parameters that produce $E_3$, with the change $\mu=0.8$.

For $P_4$, the Jacobian has the form
$$
J=
\left(
\begin{array}{ccc}
 \bigskip 
0 & - & 0\\
 \bigskip 
+ & - & - \\
? & + & +
\end{array}
\right)
\equiv
\left(
\begin{array}{ccc}
0 & Z & 0\\
B & C & D\\
E & F & G
\end{array}
\right)
$$
with
%$Z<0$, $C<0$, $D<0$, $B>0$, $F>0$, $G>0$, and 
$E$ undecided in sign. 
The characteristic equation is the cubic $\sum_{k=0}^3 c_k \lambda ^k$ with $c_3=1$, $c_2=-(C+G)$, $c_1=CG-BZ-FD$, $c_0=Z(BG-ED)$.
%\clearpage

Numerical simulations reveal that this equilibrium can be obtianed for the parameter values
$\sigma=0.4$, $r=0.6$, $\mu=0.17$, $q=0.5$, $w=0.5$, $K=15$, $m=0.33$, $g=0.14$, $f=0.2$, Figure \ref{fig_P4oscill} left.
The same parameter choice as for $E_4$, namely
$\sigma=0.5$, $r=0.5$, $\mu=0.4$, $q=0.2$, $w=0.5$, $K=3(m/g)^2$, $m=0.2$, $g=0.1$, $f=0.3$
gives instead limit cycles, Figure \ref{fig_P4oscill} right.

\begin{figure}[h]
\centering
\includegraphics[scale=0.47]{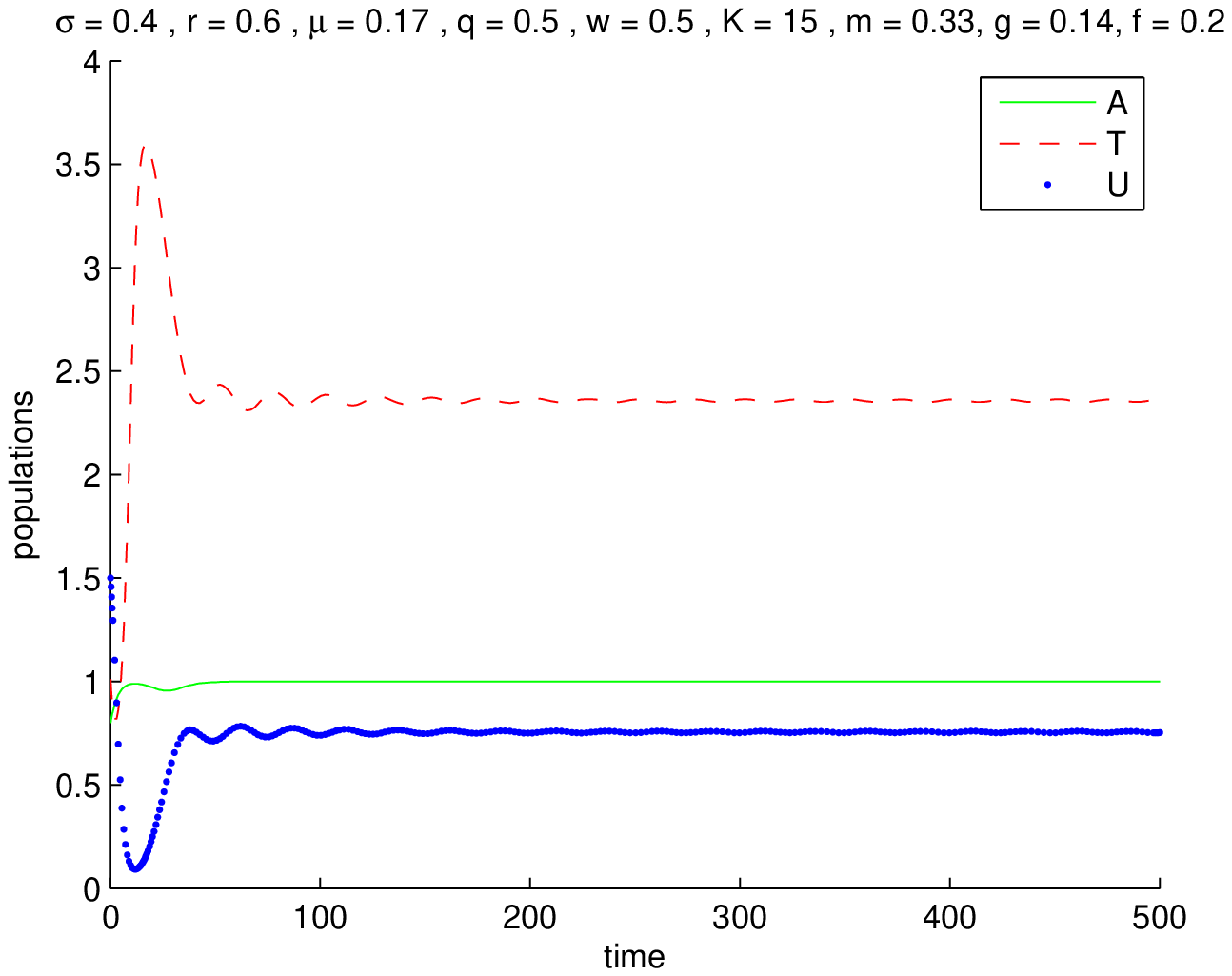}
\includegraphics[scale=0.47]{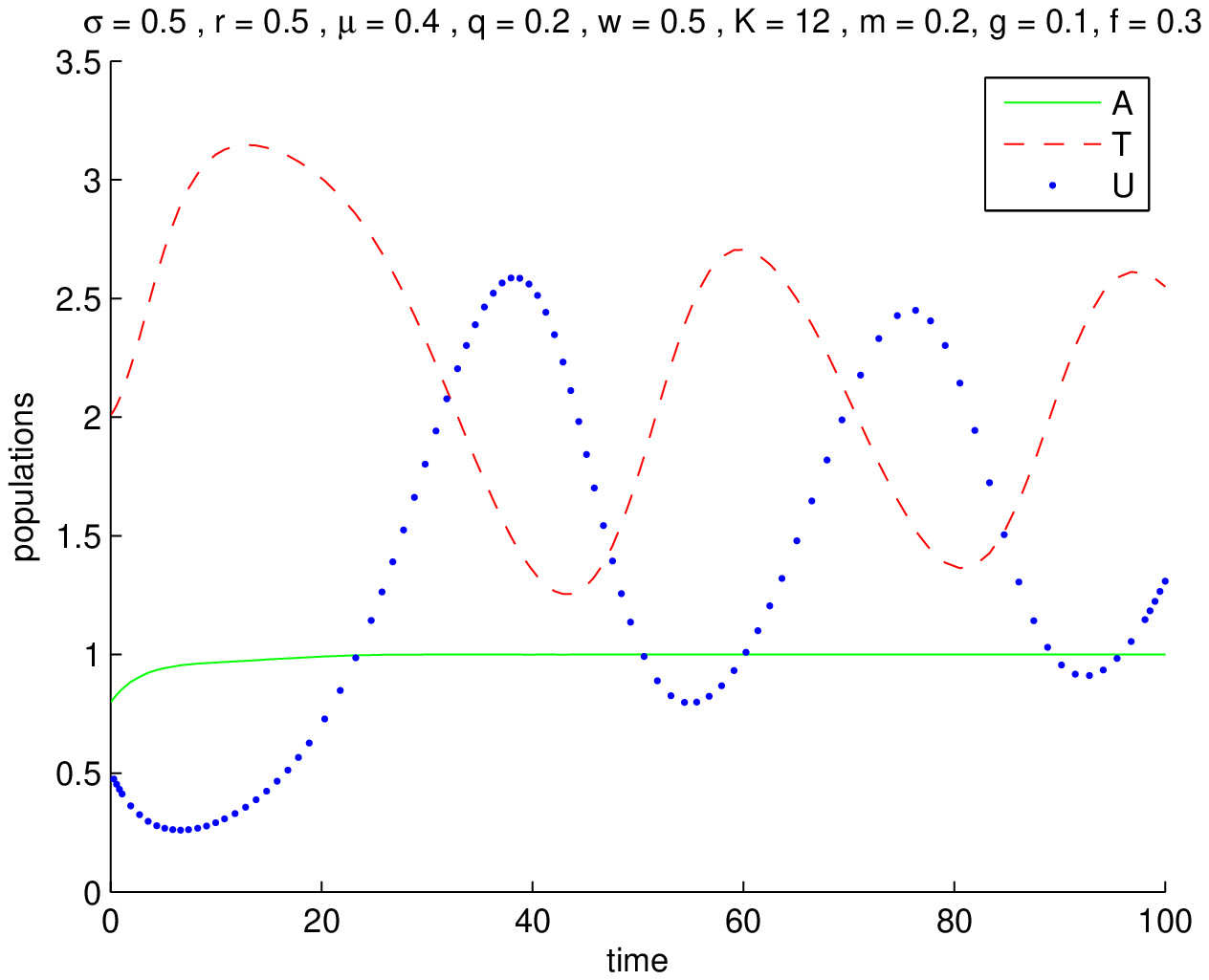}
\caption{Left: Stable equilibrium $P_4$ for the choice $\sigma=0.4$, $r=0.6$, $\mu=0.17$, $q=0.5$, $w=0.5$, $K=15$, $m=0.33$, $g=0.14$, $f=0.2$;
Right: two-dimensional limit cycle at $P_4$ obtained for the parameter values
$\sigma=0.5$, $r=0.5$, $\mu=0.4$, $q=0.2$, $w=0.5$, $K=3(m/g)^2$, $m=0.2$, $g=0.1$, $f=0.3$.}
\label{fig_P4oscill}
\end{figure}

\section{Conclusions}

In this paper we studied a predator-prey ecoepidemic model, in which an unrecoverable disease spreads by contact among the prey.
Predators feed on healthy as well as infected prey. The specific feature of this model, with respect to most of the current
literature in ecoepidemics, is that the prey gather together for defensive purposes. This herd behavior has already been introduced
in earlier investigations on demographic ecosystems, \cite{AV,APV} as well as in the case of ecoepidemics, \cite{EV11}.
But in the latter case the infected are assumed to be left behind by the herd, and hunted individually by the predators.
Here we modeled instead the case in which the diseased individuals remain in the herd, and therefore they are protected by
the ``shelter'' built by the set of prey, that gathering together in large numbers may confuse the predators.
This phenomenon has been observed in several different situations \cite{Allen,Miller,Mo,Ono,PP,TP}.
The infected are thus hunted like all the other individuals present in the herd.
Thus, mainly the individuals on the boundary of the herd suffer from the attacks of the predators.

Two cases are then examined. The first one assumes that infected prey are harmless when ingested by the predators,
while in the second one the diseased individuals are toxic for the hunting population.
To avoid mathematical difficulties due to the possible presence of singularities in the Jacobian of the system
when the prey population vanishes,
we introduced new variables in place of the original populations.
It is possible to give a very specific meaning to each one of the new depenent quantities, namely:
the ratio of healthy prey over the total amount of prey, the number of predators per prey located at the
edge of the herd area and finally the number of prey occupying the boundary of the herd.

There are only four possibly stable equilibria. In the first one, just the healthy prey thrive. Here the disease is not present, but also
the predators are wiped out.
Alternatively, while the predators still disappear, the disease remains endemic with only the prey population surviving.
Thirdly healthy prey coexist with the predators, either stably or through possibly
persistent oscillations; in this situation the disease is eradicated while the ecosystem is preserved.
The final coexistence equilibrium is also possible, with both populations thriving and an endemic disease.

Note that
the only impossible alternative in these circumstances is the predators thriving only with infected prey.
But this is a consequence of the model assumptions, because infected prey are assumed to be too weak to sustain themselves.
Their disappearance of the healthy prey prevents therefore the infected prey to replenish their population, and therefore
the latter is bound to vanish. Since the predators are specialist and do not have any source of food left, either sound or infected,
they are bound to be wiped out from the ecosystem too.

This analysis qualitatively holds for the model in which the infected are toxic for the predators.
Quantitative statements could be provided, but they should be related to specific ecosystems for which
at least some parameter values can be obtained by field measurements.

\end{document}